\documentclass[final]{amsart}
\usepackage{pdfsync}
\usepackage{amsmath,amssymb}
\usepackage{enumerate}
\usepackage{xspace}
\usepackage{boxedminipage}
\usepackage{algpseudocode}
\usepackage{listings}
\usepackage{tabularx}
\usepackage[font=footnotesize]{subfig}
\usepackage{tristan}

\renewcommand{\bih}[2]{\ensuremath{\cA_h\qp{#1,#2}}}

\renewcommand{\vec}[1]{\geovec{#1}}
\renewcommand{\mat}[1]{\geomat{#1}}

\renewcommand{\B}{\mat{B}}

\numberwithin{equation}{section}
\setboolean{showtodo}{true}
\setboolean{shownotes}{true}
\setboolean{showchanges}{false}
\setboolean{usemathrsfs}{false}
\author{
  Andreas Dedner
}
\address{
  Mathematics institute,
  University of Warwick,
  Coventry, England UK,
  CV4 7AL.}
\email{{\tt\linkedurl{A.S.Dedner@warwick.ac.uk}}.
}

\author{
  Tristan Pryer
}
\address{
  School of Mathematics, Statistics \& Actuarial Science,
  University of Kent,
  Canterbury, England UK,
  CT2 7NF.
}
\email{
    {\tt\linkedurl{T.Pryer@kent.ac.uk}}.
}

\title[DG methods for nonvariational problems] 
{Discontinuous Galerkin methods for nonvariational problems }
\date{\today}
\pdfformat{true}

\begin{document}
\maketitle
\begin{abstract}
  We extend the finite element method introduced by Lakkis and Pryer
  [2011] to approximate the solution of second order elliptic problems
  in nonvariational form to incorporate the discontinuous Galerkin (DG)
  framework. This is done by viewing the NVFEM as a mixed method
  whereby the ``finite element Hessian'' is an auxiliary variable in
  the formulation. Representing the finite element Hessian in a
  discontinuous setting yields a linear system of the same size and having
  the same sparsity pattern of the compact DG methods for variational
  elliptic problems. Furthermore, the system matrix is very easy to
  assemble,
  Thus this approach greatly reduces the computational complexity of
  the discretisation compared to the continuous approach.

  We conduct a stability and consistency analysis making use of the
  unified framework set out in Arnold et. al. [2001]. We also give an
  \apriori 
  analysis of the method. The analysis applies to any consistent
  representation of the finite element Hessian, thus is applicable to
  the previous works making use of continuous Galerkin approximations.
\end{abstract}

\section{Introduction}
\label{sec:introduction}

Nonvariational partial differential equations (PDEs) are those which
are given in the form
\begin{equation}
  \label{eq:nonvariational-problem}
  -\frob{\A}{\Hess u} = f,
\end{equation}
where $\frob{\geomat X}{\geomat Y} = \trace\qp{\Transpose{\geomat X}
  \geomat Y}$ is the Frobenious inner product between matrices. If the
matrix $\A$ is differentiable then there is an equivalence between
this problem and its variational sibling
\begin{equation}
  \label{eq:divergence-advection-form}
  -\div\qp{\A \nabla u} + \D{\A}\nabla u = f,
\end{equation}
where 
\begin{equation}
  \D\A = \qp{\sum_{i=1}^d \partial_i a_{i,1}(\geovec x) , \dots , \sum_{i=1}^d \partial_i a_{i,d}(\geovec x)}.
\end{equation}
Rewriting in this form is sometimes undesirable. For example, if the
coefficient matrix $\A$ has near singular derivatives the problem will
become advection dominated and possibly unstable for conforming finite
element methods. There is a wealth of material on the treatment of
advection dominated problems \cite[c.f.]{ErnGuermond:2004,
  ElmanSilvesterWathen:2005}. If $\A$ is not differentiable then the
problem has no variational structure. In this case standard finite
element methods cannot be applied.

In a previous work \cite{LakkisPryer:2011} a finite element method for
the approximation of the nonvariational problem
(\ref{eq:nonvariational-problem}) was introduced. This involved the
introduction of a \emph{finite element Hessian} represented in the
same finite element space as the solution (modulo boundary
conditions). The applications of the discrete representation of a
Hessian of a piecewise function are becoming broader, for example, it
can be used to drive anisotropic adaptive algorithms
\cite{AgouzalVassilevski:2002,
  ValletManoleDompierreDufourGuibault:2007}, as a notion of discrete
convexity \cite{Aguilera:2008} and in the design of finite element
methods for nonlinear fourth order problems \cite{Pryer:2012a}. We are
particularly interested in nonvariational PDEs due to their relation
to general fully nonlinear PDEs
\begin{equation}
  \funq{\Hess u} = 0,
\end{equation}
which are of significant current research. There have been finite
element methods presented for this general class of problem for
example in \cite{Bohmer:2008} the author presents a $\cont{1}$ finite
element method shows stability and consistency (hence convergence) of
the scheme which requires a high degree of smoothness on the exact
solution. In \cite{FengNeilan:2009,FengNeilan:2009a} the authors give
a method in which they approximate the general second order fully
nonlinear PDE by a sequence of fourth order quasilinear PDEs. This is
reminiscent of the vanishing viscosity method introduced for
classically studying first order fully nonlinear PDEs. Efficiency of
any method used to approximate a problem such as this is key. Each of
the methods are computationally costly due to their reliance on
$\cont{1}$ finite elements \cite{Bohmer:2008, FengNeilan:2009} or
mixed methods \cite{FengNeilan:2009a}.

In \cite{BarlesSouganidis:1991} a generic framework was set up to
prove convergence of numerical approximations to the solutions of
fully nonlinear PDEs. This involved constructing monotone sequences of
approximations which are typically applied to finite difference
approximations of the nonlinear problem \cite[c.f.]{Oberman:2006}. The
assumption of consistency made in the \cite{BarlesSouganidis:1991}
framework are incompatible with finite element methods, however, an
extremely important observation made in \cite{JensenSmears:2011} is
that the consistency condition may be weakened to incorporate the
finite element case using a \emph{localisation argument} (in the case
of isotropic diffusion). 

In this contribution we present a method for the discontinuous
approximation of the linear nonvariational problem
(\ref{eq:nonvariational-problem}). We also present convergence
analysis for a certain subclass of the nonvariational problems, those
which are coercive. This allows us to use variational techniques to
analyse the problem. We prove optimal convergence rates for the finite
element solution in broken Sobolev norms. Note that the results
presented here are immediately applicable to the method derived for
the continuous case given in \cite{LakkisPryer:2011}.

The algebraic formulation of the continuous approximation of the
nonvariational problem requires the solution of large sparse
$\qp{d+1}^2\times N^2$ linear system \cite[Lem 3.3]{LakkisPryer:2011},
where $d$ is the dimension of the problem and $N$ the number of
degrees of freedom. Equivalently, using a Sch\"ur complement argument,
this can be reduced to an $N^2$ full linear system. The reason that
this system is full is due to the global nature of the $\leb{2}(\W)$
projection operator into a continuous finite element space. The
motivation for extending the nonvariational finite element method into
the discontinuous setting is the massive gain in computational
efficiency over the continuous case. Indeed, due to the local
representation of the projection operators in these discontinuous
spaces we are able to make massive computational savings, in that the
system matrix will become sparse and is the same size as that of a
standard discontinuous Galerkin stiffness matrix.


To test the method numerically we make use of the finite element
package \dune
\cite{BastianBlattDednerEngwerKlofkornKornhuberOhlbergerSander:2008,
  BastianBlattDednerEngwerKlofkornOhlbergerSander:2008}. In this work
we are interested in the asymptotic behaviour of the discontinuous
approximation. In a subsequent work we will study the computational
gains using the discontinuous framework presented over the continuous
one given in \cite{LakkisPryer:2011}, as well as exploit the powerful
parallelisation capabilities of the package.

The rest of the paper is set out as follows: In \S\ref{sec:notation}
we formally introduce the model problem and give a brief review of
known classical facts about nonvariational PDEs. In
\S\ref{sec:discretisation} we examine the discretisation of the
nonvariational method in the discontinuous Galerkin framework, making
use of the unified framework set out in
\cite{ArnoldBrezziCockburnMarini:2001} to derive a very general
formulation of the finite element Hessian represented as a
discontinuous object. We present some examples and examine the natural
question of what happens when we try to eliminate the finite element
Hessian from the formulation. In \S\ref{sec:consistency} we look at
the consistency and stability of the finite element Hessian and
present our main analytical results of convergence. Finally, in
\S\ref{sec:numerics} we detail a summary of extensive numerical
experiments aimed at examining convergence and computational speed of
the method presented.


\section{Problem formulation}
\label{sec:notation}

In this section we formulate the model problem, fix notation and give
some basic assumptions. In addition we review the existence and
uniqueness of the nonvariational problems. We begin by introducing the
Lebesgue spaces
\begin{gather}
  \leb{2}(\W) = \ensemble{\phi} {\int_\W \norm{\phi(\geovec x)}^2 \d \geovec x < \infty}
  \AND
  \leb{\infty}(\W) 
  =
  \ensemble{\phi} {\sup_{\geovec x \in \W} \norm{\phi(\geovec x)} < \infty},
\end{gather}
and the Sobolev and Hilbert spaces
\begin{gather}
  \sob{k}{p}(\W) 
  = 
  \ensemble{\phi\in\leb{p}(\W)}
  {\D^{\vec\alpha}\phi\in\leb{p}(\W), \text{ for } \norm{\geovec\alpha}\leq k} \AND
  \sobh{k}(\W)
  := 
  \sob{k}{2}(\W).
\end{gather}
These are equipped with the norms
\begin{gather}
  \Norm{\phi}^2_{\leb{2}(\W)} = \int_\W \norm{\phi}^2 \d \geovec x
  , \qquad 
  \Norm{\phi}_{\leb{\infty}(\W)} = \sup_{\geovec x \in \W} \norm{\phi(\geovec x)},
  \\
  \Norm{v}_{\sob{k}{p}(\W)}^p
  = 
  \sum_{\norm{\vec \alpha}\leq k}\Norm{\D^{\vec \alpha} v}_{\leb{p}(\W)}^p 
  \AND 
  \norm{v}_{\sob{k}{p}(\W)}^p
  =
  \sum_{\norm{\vec \alpha} = k}\Norm{\D^{\vec \alpha} v}_{\leb{p}(\W)}^p.
\end{gather}
where $\vec\alpha = \{ \alpha_1,...,\alpha_d\}$ is a multi-index,
$\norm{\vec\alpha} = \sum_{i=1}^d\alpha_i$ and derivatives
$\D^{\vec\alpha}$ are understood in a weak sense. We pay particular
attention to the cases $k = 1,2$ and
\begin{gather}
  \hoz(\W) := \text{closure of }\cont{\infty}_0(\W) \text{ in } \sobh{1}(\W).
\end{gather}

The model problem in strong form is: 
Find $u\in\sobh2(\W)\cap\hoz(\W)$ such that
\begin{equation}
  \begin{split}
    \label{Problem}
    \ltwop{\linop u}{\phi} &= \ltwop{f}{\phi} \qquad \Foreach
    \phi\in\hoz(\W), \\
  \end{split}
\end{equation}
where the data {$f\in\leb2(\W)$} is prescribed and $\linop$ is a
general linear, second order, uniformly elliptic partial differential
operator. {Let $\A \in \Le{\infty}(\W)^{d\times d}$, 
  we then define}
\begin{equation}
  \label{eqn:def:nondivergence-linear-operator}
  \dfunkmapsto[.]
	      {\linop}
	      u
	      {\sobh{2}(\W)\cap\hoz(\W)}
	      {\linop u:= -\frob{\A}{\Hess u}}
	      {\leb{2}(\W)}
\end{equation}

We assume that $\A$ is uniformly positive definite, i.e., there exists
a $\gamma>0$ such that for all $\vec x$
\begin{equation}
  \label{eq:ellipticity}
  \Transpose{\vec y} \A(\vec x) \vec y \geq \gamma \norm{\vec y}^2 
  \Foreach \vec y\in\reals^d,
\end{equation}
and we call $\gamma$ the \emph{ellipticity constant}. 


Nonvariational PDEs are not as well studied as their variational
brethren from a numerical analysis viewpoint. For the benefit of the
reader we present a concise account of known results for strong
solutions of this class of problem.

\index{strong solution}
\begin{Defn}[strong solution]
  A \emph{strong solution} of (\ref{eq:nonvariational-problem}) is a
  function $u\in\sobh2(\W) \cap \hoz(\W)$, that is a twice weakly
  differentiable function, which satisfies the problem almost
  everywhere.
\end{Defn}
\begin{The}[existence and regularity of a strong solution of (\ref{eq:nonvariational-problem}) {\cite[Thm 9.15]{Gilbarg:1983}}]
  \label{the:regularity-of-nonvar-form}
  Let $\W\subset\reals^d$ be a $\cont{1,1}$ domain. Suppose also that
  $\A\in\cont{0}(\W)^{d\times d}$ and $f\in\leb{2}(\W)$ such that the
  problem
  \begin{equation}
    \label{eq:existence-of-strong-soln}
    \begin{split}
      -\frob{\A}{\Hess u} &= f \text{ in }\W
      \\
      u&=0 \text{ on }\partial \W
    \end{split}
  \end{equation}
  {is uniformly elliptic}. Then (\ref{eq:existence-of-strong-soln})
  has a unique strong solution. There also exists a constant
  independent of $u$ such that
  \begin{equation}
    \Norm{u}_{\sobh{2}(\W)}\leq C\Norm{f}_{\leb{2}(\W)}.
  \end{equation}
\end{The}

\begin{Rem}[less regular solutions]
  Note that the theory of viscosity solutions has been developed for
  non classical solutions of (\ref{Problem}) if the problem data does
  not satisfy the regularity assumed above see \cite{Gilbarg:1983}.
\end{Rem}

\begin{Hyp}[regularity of $\A$]
  From hereon in we will assume that the problem data is sufficiently
  smooth such that solutions exist and belong to at least
  $\sobh{2}(\W)\cap\hoz(\W)$.
\end{Hyp}
\begin{Rem}[regularity of $\W$]
  Theorem \ref{the:regularity-of-nonvar-form} specifies that $\W$ must
  be a $\cont{1,1}$ domain. We will be approximating such a domain
  with one which is only $\cont{0,1}$ (\ie a polyhedral one). We thus
  assume that the model problem admits a unique strong solution even
  when $\W$ is only $\cont{0,1}$. To circumvent this assumption curved
  finite elements could be used to fit the boundary exactly
  \cite{bernardi:1989}. For simplicity we will not present this case
  here, although we believe our analysis can be extended to this case.
\end{Rem}

\section{Discretisation}
\label{sec:discretisation}

Let $\T{}$ be a conforming, shape regular triangulation of $\W$,
namely, $\T{}$ is a finite family of sets such that
\begin{enumerate}
\item $K\in\T{}$ implies $K$ is an open simplex (segment for $d=1$,
  triangle for $d=2$, tetrahedron for $d=3$),
\item for any $K,J\in\T{}$ we have that $\closure K\meet\closure J$ is
  a full subsimplex (i.e., it is either $\emptyset$, a vertex, an
  edge, a face, or the whole of $\closure K$ and $\closure J$) of both
  $\closure K$ and $\closure J$ and
\item $\union{K\in\T{}}\closure K=\closure\W$.
\end{enumerate}
We use the convention where $\funk h\W\reals$ denotes the
\emph{meshsize function} of $\T{}$, i.e.,
\begin{equation}
  h(\vec{x}):=\max_{\closure K\ni \vec x}h_K,
\end{equation}
where $h_K$ is the diameter of $K$. We let $\E{}$ be the skeleton (set
of common interfaces) of the triangulation $\T{}$ and say $e\in\E$ if
$e$ is on the interior of $\W$ and $e\in\partial\W$ if $e$ lies on the
boundary $\partial\W$.
The assumptions on the tessellation made here are typical in the finite
element analysis. For the presentation of the method and its analysis, 
some assumption could be relaxed 
(e.g. the form of the elements or the assumption on a conforming
triangulation) but this would lead to an unnecessary increase in the complexity 
of the presentation.

Let $\poly k(\T{})$ denote the space of piecewise polynomials of
degree $k$ over the triangulation $\T{}$,\ie
\begin{equation}
  \poly k (\T{}) = \ensemble{ \phi }{\phi|_K \in \poly k (K) }
\end{equation}
 and introduce the \emph{finite element spaces}
\begin{gather}
  \label{eqn:def:finite-element-space}
  \dgzero = \dgzero\qp{\T{},k} := \poly k(\T{}) \cap \hoz(\W)
  \\
  \dg = \dg\qp{\T{},k} := \poly k(\T{}) 
\end{gather}
to be the usual spaces of discontinuous piecewise polynomial
functions which are compactly and non compactly supported respectively.

\begin{Rem}[generalised Hessian]
  \label{rem:generalised-hessian}
  Assume a function $v\in\sobh2(\W)$, let $\geovec
  n:\partial\W\to\reals^d$ be the outward pointing normal of $\W$ then
  the Hessian $\Hess v$ of $v$, satisfies the following identity:
  \begin{equation}
    \label{eq:generalised-hessian}
    \int_\W{\Hess v} \ {\phi} \d \geovec x
    = 
    -
    \int_\W{\nabla v}\otimes{\nabla \phi} \d \geovec x
    +
    \int_{\partial\W}{\nabla v}\otimes{\geovec n \ \phi} \d s \Foreach \phi\in\sobh1(\W).
  \end{equation}
  If $v\in\sobh{1}(\W)$ (\ref{eq:generalised-hessian}) is still well defined in view of
  duality, in this case we set
  \begin{equation}
    \duality{\Hess v}{\phi} 
    =
    -
    \int_\W{\nabla v}\otimes{\nabla \phi}\d \geovec x
    +
    \int_{\partial\W}{\nabla v}\otimes{\geovec n \ \phi} \d  s \Foreach \phi\in\sobh1(\W),
  \end{equation}
  where the last term is understood as a pairing between
  $\sobh{-1/2}(\W)$ and $\sobh{1/2}(\W)$.
\end{Rem}

\begin{Defn}[broken Sobolev spaces, trace spaces]
  \label{defn:broken-sobolev-space}
  We introduce the broken Sobolev space
  \begin{equation}
    \sobh{k}(\T{})
    :=
    \ensemble{\phi}
             {\phi|_K\in\sobh{k}(K), \text{ for each } K \in \T{}}.
  \end{equation}
  We also make use of functions defined in these broken spaces
  restricted to the skeleton of the triangulation. This requires an
  appropriate trace space
  \begin{equation}
    \Tr{\E} := \prod_{K\in\T{}} \leb{2}(\partial K) = \prod_{K\in\T{}} \sobh{\frac{1}{2}}(K).
  \end{equation}
\end{Defn}

\begin{Defn}[jumps, averages and tensor jumps]
  \label{defn:averages-and-jumps}
  We define average, jump and tensor jump operators for arbitrary
  scalar functions $v\in\Tr{\E}$, vectors $\vec
  v\in\Tr{\E}^d$ and matrices $\mat
  V\in\Tr{\E}^{d\times d}$ as
  \begin{gather}
    \label{eqn:average}
    \avg{v} =  {\frac{1}{2}\qp{v|_{K_1} + v|_{K_2}}},
    \qquad \avg{\vec v} = {\frac{1}{2}\qp{\vec{v}|_{K_1} + \vec{v}|_{K_2}}},
    \\\nonumber\\
    \label{eqn:jump}
    \jump{v} = {{{v}|_{K_1} \geovec n_{K_1} + {v}|_{K_2}} \geovec n_{K_2}},
    \qquad
    \jump{\vec v}
    = 
    \Transpose{\qp{\vec{v}|_{K_1}}}\geovec n_{K_1} + \Transpose{\qp{\vec{v}|_{K_2}}}\geovec n_{K_2},
    \\\nonumber\\
    \label{eqn:jump-mat}
    \jump{\mat V}
    =
    {{{{\mat{V}|_{K_1}}}\geovec n_{K_1} + {{\mat{V}|_{K_2}}}\geovec n_{K_2}}},
    \qquad
    \tjump{\vec v}
    =
    {{\vec{v}|_{K_1} }\otimes \geovec n_{K_1} + \vec{v}|_{K_2} \otimes\geovec n_{K_2}}.
  \end{gather}

  Note that on the boundary of the domain $\partial\W$ the jump and
  average operators are defined as
  \begin{gather}
    \avg{v}
    \Big\vert_{\partial\W} 
    := v,
    \qquad 
    \avg{\geovec v}
    \Big\vert_{\partial\W}
    :=
    \geovec v,
    \\
    \jump{v}
    \Big\vert_{\partial\W}
    := 
    v\geovec n,
    \qquad 
    \jump{\geovec v}
    \Big\vert_{\partial\W} 
    :=
    \Transpose{\geovec v}\geovec n,
    \\
    \jump{\mat V}
    \Big\vert_{\partial\W} 
    :=
    \mat V \geovec n,
    \qquad
    \tjump{\vec v}
    \Big\vert_{\partial\W} 
    :=
    \vec v\otimes \geovec n.
  \end{gather}
\end{Defn}

We will often use the following Proposition which we state in full for
clarity but whose proof is merely using the identities in Definition
\ref{defn:averages-and-jumps}.
\begin{Pro}[elementwise integration]
  \label{Pro:trace-jump-avg}
  For a generic vector valued function $\geovec p$ and scalar valued
  function $\phi$ we have
  \begin{equation}
    \label{eq:jump-avg-eq1}
    \begin{split}
      \sum_{K\in\T{}}
      \int_K \div\qp{\geovec p} \phi \d \geovec x
      =
      \sum_{K\in\T{}}
      \qp{
        -
        \int_K
        \Transpose{\geovec p} \nabla_h \phi \d \geovec x
        +
        \int_{\partial K}
        \phi \Transpose{\geovec p} \geovec n_K \d s
        },
    \end{split}
  \end{equation}
  where $\nabla_h = \Transpose{\qp{\D_h }}$ is the elementwise spatial
  gradient. Furthermore, If we have $\geovec p \in \Tr{\E\cup\partial\W}^d$ and $\phi \in
  \Tr{\E\cup\partial\W}$, the following identity holds
  \begin{equation}
    \label{eq:jump-avg-eq2}
    \sum_{K\in\T{}}
    \int_{\partial K}
    \phi \Transpose{{\geovec p}}
    \geovec n_K \d s
    =
    \int_\E 
    \jump{\geovec p} 
    \avg{\phi}
    \d s
    +
    \int_{\E\cup\partial\W}
    \Transpose{\jump{\phi}}
    \avg{\geovec p}
    \d s
    =
    \int_{\E\cup\partial\W}
    \jump{\geovec p \phi}
    \d s,
  \end{equation}
  An equivalent tensor formulation of
  (\ref{eq:jump-avg-eq1})--(\ref{eq:jump-avg-eq2}) is
  \begin{equation}
    \label{eq:jump-avg-eq3}
    \begin{split}
      \sum_{K\in\T{}}
      \int_K \D_h {\geovec p} \phi \d \geovec x
      =
      \sum_{K\in\T{}}
      \qp{
        -
        \int_K
        {\geovec p} \otimes \nabla_h \phi \d \geovec x
        +
        \int_{\partial K}
        \phi {\geovec p} \otimes \geovec n_K \d s
        },
    \end{split}
  \end{equation}
  where
  \begin{equation}
    \label{eq:jump-avg-eq4}
    \sum_{K\in\T{}}
    \int_{\partial K}
    \phi {\geovec p}
    \otimes
    \geovec n_K \d s
    =
    \int_\E 
    \tjump{\geovec p} 
    \avg{\phi}
    \d s
    +
    \int_{\E\cup\partial\W}
    {\jump{\phi}}\otimes
    \avg{\geovec p}
    \d s
    =
    \int_{\E\cup\partial\W}
    \tjump{\geovec p \phi}
    \d s.
  \end{equation}
  In addition for matrix valued $\mat V$ we have that
  \begin{equation}
    \sum_{K\in\T{}}
    \int_{K}
    \frob{
      \qp{\D_h {\geovec p}}
    }
     {\mat V} \d \vec x
     =
     \sum_{K\in\T{}}
     \qp{ -\int_{K}
     \frob{
       {{\geovec p}}
    }
     {\D_h \mat V}\d \vec x
     +
     \int_{\partial \W}
     \Transpose{\qp{\mat V \vec p}}\vec n \d s
   }
   \end{equation}
   and
  \begin{equation}
    \sum_{K\in\T{}}
    \int_{\partial \W}
     \Transpose{\qp{\mat V \vec p}}\vec n
     \d s
    =
    \int_\E 
    \Transpose{\jump{\mat V}}\avg{\vec p} 
    \d s
    +
    \int_{\E\cup\partial\W}
    \frob{\tjump{\vec p}}{\avg{\mat V}}
    \d s
    =
    \int_{\E\cup\partial\W}
    \jump{\mat V\geovec p}
    \d s.
  \end{equation}  
\end{Pro}

\subsection{Construction of an appropriate discrete Hessian}

We now use the framework set out in
\cite{ArnoldBrezziCockburnMarini:2001} to construct a general notion
of discrete Hessian. We first give a definition using a flux formulation:
\begin{Defn}[generalised finite element Hessian: flux formulation]
  \label{the:fully-generalised-fe-hessian-flux-form}
  Let $u\in\sobh2(\T{})$ and
  $\hat U : \sobh1(\T{}) \to \Tr{\E\cup\partial\W}$ be a linear form and
  $\hat{\geovec p} : \sobh2(\T{}) \times \sobh1(\T{})^d \to
  \Tr{\E\cup\partial\W}^d$ a bilinear form representing approximations
  to $u$ and $\nabla u$ over the skeleton of the triangulation.
  Then we define the generalized finite element Hessian
  $\H[u]$ as the solution of
  \begin{gather}
    \label{eq:primal-hessian}
    \int_K {\H[u]} \ {\Phi} \d \geovec x
    =
    -
    \int_K \geovec p \otimes \nabla_h \Phi\d \geovec x
    +
    \int_{\partial K} \geovec{\hat{p}}_K \otimes \geovec n \ \Phi\d s
    \Foreach \Phi \in \sobh{1}(\T{})
    \\
    \label{eq:primal-hessian-2}
    \int_K \geovec p \otimes \geovec q \d \geovec x
    =
    -
    \int_K u \ \D_h \geovec q \d \geovec x
    +
    \int_{\partial K} \geovec q \otimes \geovec n  \ \hat{U}_K \d s
    \Foreach \geovec q \in \qp{\sobh{1}(\T{})}^d,
  \end{gather}
  for all $\Phi\in\dg$.
\end{Defn}

We now present the primal formulation for the generalized finite element
Hessian:

\begin{The}[generalised finite element Hessian: primal form]
  \label{the:fully-generalised-fe-hessian}
  Let $u\in\sobh2(\T{})$ and let $\hat U$ and $\hat{\geovec p}$ be
  defined as in Definition
  \ref{the:fully-generalised-fe-hessian-flux-form}, then the
  generalised finite element Hessian $H[u]$ is given for each
  $\Phi\in\dg$ as
  \begin{equation}
    \label{eq:fully-generalised-fe-hessian}
    \begin{split}
      \int_{\W} \H[u] \ \Phi\d \geovec x
      &=
      -
      \int_\W \nabla_h u \otimes \nabla_h \Phi\d \geovec x
      +
      \int_{\E\cup \partial\W} \jump{\Phi}\otimes\avg{\geovec{\hat p}}\d s
      +
      \int_\E \avg{\Phi} \tjump{\geovec{\hat p}}\d s
      \\
      &\qquad- 
      \int_{\E} \avg{\hat U - u} \tjump{\nabla_h \Phi}\d s
      -
      \int_{\E \cup \partial \W} \jump{\hat U - u} \otimes \avg{\nabla_h \Phi}\d s.
    \end{split}
  \end{equation} 
\end{The}

\begin{Proof}
Note that in view of Definition \ref{defn:averages-and-jumps} for
generic vector fields $\geovec q \in \fesW$ and $v \in \fes$ we have
the following identity
\begin{equation}
  \label{eq:edge-identity}
  \sum_{K\in\T{}} \int_{\partial K} v \geovec q \otimes \geovec n \d s
  =
  \int_{\E \cup \partial \W} \jump{v}\otimes \avg{\geovec q} \d s
  +
  \int_\E \avg{v} \tjump{\geovec q}\d s.
\end{equation}
Then summing (\ref{eq:primal-hessian}) over $K\in\T{}$ and making use
of the identity (\ref{eq:edge-identity}) we see
\begin{equation}
  \begin{split}
    \int_\W {\H[u]} \ {\Phi} \d \geovec x
    &=
    \sum_{K\in\T{}}
    \int_K {\H[u]} \ {\Phi} \d \geovec x
    =
    \sum_{K\in\T{}}   
    \qp{-  
    \int_K \geovec p \otimes \nabla_h \Phi\d \geovec x
    +
    \int_{\partial K} \geovec{\hat{p}}_K \otimes \geovec n \ \Phi}\d s
    \\
    &=
    -  
    \int_\W \geovec p \otimes \nabla_h \Phi\d \geovec x
    +
    \int_{\E \cup \partial \W} \jump{\Phi} \otimes \avg{\geovec{\hat{p}}_K} \d s
    +
    \int_\E \avg{\Phi} \tjump{\geovec{\hat{p}}_K}\d s.
  \end{split}
\end{equation}
Using the same argument for (\ref{eq:primal-hessian-2})
\begin{equation}
  \begin{split}
    \int_\W \geovec p \otimes \geovec q \d \geovec x
    &=
    \sum_{K\in\T{}}
    \int_K \geovec p \otimes \geovec q \d \geovec x
    =
    \sum_{K\in\T{}}
    \qp{
      -
      \int_K u \ \D_h \geovec q \d \geovec x
      +
      \int_{\partial K} \geovec q \otimes \geovec n  \ \hat{U}_K \d s
    }
    \\
    &=
    -
    \int_\W u \ \D_h \geovec q \d \geovec x
    +
    \int_{\E \cup \partial \W} \jump{\hat{U}} \otimes \avg{\geovec q}\d s
    +
    \int_\E \avg{\hat{U}} \tjump{\geovec q}\d s.
  \end{split}
\end{equation}

Note that, again making use of (\ref{eq:edge-identity}) we have for
each $\geovec q\in\sobh1(\T{})^d$ and $v\in \sobh1(\T{})$ that
\begin{equation}
  \label{eq:edge-identity-2}
  \int_\W \geovec q \otimes \nabla_h v \d \geovec x
  =
  -
  \int_\W \D_h \geovec q v\d \geovec x
  +
  \int_{\E\cup \partial \W} \avg{\geovec q}\otimes \jump{v}\d s
  +
  \int_\E \tjump{\geovec q} \avg{v}\d s.
\end{equation}
Taking $v=u$ in (\ref{eq:edge-identity-2}) and substituting into
(\ref{eq:primal-hessian-2}) we see
\begin{equation}
  \label{eq:primal-hessian-3}
  \int_\W \geovec p\otimes \geovec q\d \geovec x
  =
  \int_\W \geovec q \otimes \nabla_h u\d \geovec x
  +
  \int_{\E \cup \partial \W} \jump{\hat U - u}\otimes \avg{\geovec q}\d s
  +
  \int_\E \avg{\hat U - u} \tjump{\geovec q}\d s.
\end{equation}
Now choosing $\geovec q = \nabla_h \Phi$ and substituting
(\ref{eq:primal-hessian-3}) into (\ref{eq:primal-hessian}) we arrive
at the fully generalised finite element Hessian given by
(\ref{eq:fully-generalised-fe-hessian}).
\end{Proof}

\begin{Rem}[consistent representations of the gradient operator]
  If one were interested in consistent representations of other
  derivatives, for example the gradient operator, one would need to
  modify the proof of Theorem
  \ref{the:fully-generalised-fe-hessian}. Examples of consistent
  gradient representations can be found in
  \cite{ArnoldBrezziCockburnMarini:2001}. See also
  \cite{BuffaOrtner:2009, Di-PietroErn:2010, BurmanErn:2008}. Using
  this methodology it should be possible to construct an entire
  hierarchy of derivatives.
\end{Rem}
\begin{Example}
  \label{ex:ip-hessian}
  An example of a dG formulation for the approximation to the Hessian,
  $\Hess u$, can be derived by taking the fluxes in the following way
  \begin{gather}
    \hat U = \begin{cases}
      \theta\avg{u_h} \text{ over } \E
      \\
      0 \text{ on } \partial \W
      \end{cases}
    \\
    \hat{\geovec p}
    =
    \avg{\nabla_h u_h} \text{ on } \E\cup \partial\W,
  \end{gather}
  where $\theta\in\{-1,1\}$. The result is 
  a discrete representation of the Hessian $\H[u_h]$ as unique element
  of $\dg^{d\times d}$ such that
  \begin{equation}
    \begin{split}
      \int_{\W} \H[u_h] \ \Phi\d \geovec x
      &=
      -
      \int_\W \nabla_h u_h \otimes \nabla_h \Phi\d \geovec x
      \\
      &\qquad +
      \int_{\E \cup \partial \W} 
      \theta \jump{u_h} \otimes \avg{\nabla_h \Phi}
      +
      \jump{\Phi}\otimes\avg{\nabla_h u_h}\d s
      \\
      &=
      \int_\W \Hess_h u_h \Phi\d \geovec x
      -
      \int_\E \tjump{\nabla_h u_h} \avg{\Phi}\d s
      \\
      &\qquad +
      \int_{\E\cup\partial\W} \theta \jump{u_h} \otimes \avg{\nabla_h \Phi} \d s \Foreach \Phi\in\dg.
    \end{split}
  \end{equation} 
\end{Example}

\subsection{The discontinuous nonvariational finite element method}
We are now in a position to state the numerical method for the
approximation of (\ref{eq:nonvariational-problem}). We look to find
$u_h\in\dgzero$ together with $\H[u_h]\in\dg^{d\times d}$ such that
\begin{gather}
  \label{eq:nonvariational-fe-approx}
  \bih{u_h}{\Psi} = l(\Psi) \Foreach \Psi\in\dgzero
\end{gather}
with
\begin{gather}
  \bih{u_h}{\Psi} := 
  \int_\W -\frob{\A}{\H[u_h]} \Psi \d \geovec x
  +
  \int_{\E\cup\partial\W} \sigma h^{-1} \Transpose{\jump{u_h}}\jump{\Psi}\d s
  \\
  l(\Psi) := \int_\W f \Psi \d \geovec x,
\end{gather}
where the \emph{penalisation parameter} $\sigma>0$ is to be chosen
sufficiently large to guarantee coercivity.

Using the $\leb{2}$ projection operator $\ltwoproj:\leb{2}(\W) \to \dg$
defined for $v\in \leb{2}(\W)$ through
\begin{equation}
\int_\W \ltwoproj\qp{v}\Psi\d \geovec x = \int_\W v\Psi \d \geovec x\quad\forall\Psi\in\dg\label{eq:l2-proj-def}
\end{equation}
it is possible to elliminate the finite element Hessian from the bilinear
form for sufficiently smooth $\A$:
\begin{Lem}[elimination of the finite element Hessian in a general setting]
  \label{the:elem-of-fe-hess}
  If $\A\in\qb{\sob{k+1}{\infty}(\W)}^{d\times d}$  and the fluxes are chosen as in
  Example~\ref{ex:ip-hessian} then
  \begin{equation}
    \begin{split}
      \bih{u_h}{\Psi} &= 
      \int_\W D_h\qp{\ltwoproj\qp{\Psi \A}}\nabla_h u_h\d \geovec x
      -
      \int_{\E\cup\partial\W} \theta \Transpose{\jump{u_h}} \avg{D_h\qp{\ltwoproj\qp{\Psi \A}}}\d s
      \\
      &\qquad -
      \int_{\E\cup\partial\W} \Transpose{\jump{\ltwoproj\qp{\Psi \A}}} \avg{\nabla_h u_h} \d s
      +
      \int_{\E\cup\partial\W} \sigma h^{-1} \Transpose{\jump{u_h}}\jump{\Psi}\d s.
    \end{split}
  \end{equation}
\end{Lem}
\begin{proof}
  This follows from the following identity
  \begin{equation}
    \begin{split}
      \int_\W -\frob{\A}{\H[u_h]} \Psi\d \geovec x
      &= 
      \int_\W -\frob{\H[u_h]}{\qp{\Psi\A}} \d \geovec x
      = 
      \int_\W -\frob{\H[u_h]}{\ltwoproj\qp{\Psi\A}} \d \geovec x
      \\
      &=
      \int_\W D_h\qp{\ltwoproj\qp{\Psi \A}}\nabla_h u_h\d \geovec x
      -
      \int_{\E\cup\partial\W} \theta \Transpose{\jump{u_h}} \avg{D_h\qp{\ltwoproj\qp{\Psi \A}}}\d s
      \\
      &\qquad -
      \int_{\E\cup\partial\W} \Transpose{\jump{\ltwoproj\qp{\Psi \A}}} \avg{\nabla_h u_h}\d s.
    \end{split}
  \end{equation}
\end{proof}  
\begin{Rem}
  The solution of the problem in this form is nontrivial due to the
  global $\leb{2}(\W)$ projection appearing in the
  formulation. However, in the discontinuous setting the global
  $\leb{2}(\W)$ projection is in fact computable locally. We may
  actually exploit this fact to optimise our schemes efficiency. We
  will discuss this further in the sequel.
\end{Rem}


\begin{Example}[Laplacian formulation]
  \label{eq:laplace-formulation}
  Note that if in (\ref{eq:nonvariational-problem}) we have that $\A =
  \geovec I$ then we have that
  \begin{equation}
    f = -\frob{A}{\Hess u} = -\Delta u
  \end{equation}
  and our bilinear form reduces to
  \begin{equation}
    \begin{split}
      \bih{u_h}{\Psi} &= 
      \int_\W \Transpose{\qp{\nabla_h\Psi}} \nabla_h u_h\d \geovec x
      -
      \int_{\E\cup\partial\W} \theta \Transpose{\jump{u_h}} \avg{\nabla_h\Psi }\d s
      \\
      &\qquad -
      \int_{\E\cup\partial\W} \Transpose{\jump{\Psi}\avg{\nabla_h u_h}} -
      \sigma h^{-1} \Transpose{\jump{u_h}}\jump{\Psi} \d s
    \end{split}
  \end{equation}
  since $\ltwoproj\qp{\Psi \A}=\Psi \geomat I$.

  The nonvariational finite element method thus coincides with the
  classical (symmetric) interior penalty method for the Laplacian
  \cite{DouglasDupont:1976}.
\end{Example}

\begin{Rem}[relation to standard dG methods]
  It is not difficult to prove that choosing to numerical fluxes in
  the same way as presented in \cite[Table
  3.2]{ArnoldBrezziCockburnMarini:2001} results in the same
  correlation to the dG methods summarised in the aforementioned paper
  for the case that $\A$ is constant. For brevity we will
  not prove this here.

  Note that when $\A$ is not constant we have that the
  nonvariational finite element method does \emph{not} coincide with
  its standard variational finite element counterpart. There is an
  extra stability property which allows the method to successfully
  cope with advection dominated problems
  \cite[\S4.2]{LakkisPryer:2011} which is illustrated by the result of Lemma
  \ref{the:elem-of-fe-hess}.
\end{Rem}

We conclude this section with a proof consistency of the method and 
then show that Galerkin orthogonality holds.
\begin{Lem}[consistency]
  Let $u\in\sobh2(\T{})$ and assume that the numerical fluxes are chosen
  in a consistent fashion in the sense of
  \cite[\S 3.1]{ArnoldBrezziCockburnMarini:2001}, that is,
  \begin{gather}
    \hat U = u|_{\E \cup \partial\W}
    \\
    \hat{\geovec p} = \nabla u|_{\E \cup \partial\W}
  \end{gather}
  Then for $\Phi\in\dg$
  \begin{equation}
    \begin{split}
      \int_\W \H[u] \Phi \d \geovec x&= \int_\W \Hess u \Phi \d \geovec x
    \end{split}
  \end{equation}
  Therefore we have that $\H[u] = \ltwoproj\qp{\Hess u}$.
\end{Lem}
\begin{Proof}
  Applying Proposition \ref{Pro:trace-jump-avg} to the first term in
  the definition of $\H[u]$ yields
  \begin{equation}
    \begin{split}
      \int_\W \H[u] \Phi \d \geovec x
      &= 
      \int_\W \Hess u \Phi \d \geovec x
      +
      \int_\E
      \tjump{\hat{\vec p}-\nabla u}
      \avg{\Phi} \d s
      +
      \int_{\E\cup \partial\W}
      \avg{\geovec{\hat p}-\nabla u}\otimes\jump{\Phi}\d s
      \\
      &\qquad
      -
      \int_{\E} \avg{\hat U - u}
      \tjump{\nabla_h \Phi}\d s
      -
      \int_{\E \cup \partial \W} \jump{\hat U - u}
      \otimes
      \avg{\nabla_h \Phi} \d s
      \\
      &= \int_\W \Hess u \Phi \d \geovec x\Foreach \Phi\in\dg.
    \end{split}
  \end{equation}
  which proves the results under the consistency conditions on the fluxes.
\end{Proof}
\begin{Lem}[Galerkin orthogonality]
  \label{lem:gal-orthog}
  Let $u\in\sobh{2}(\W)\cap\hoz(\W)$ be a strong solution to the
  problem (\ref{eq:nonvariational-problem}) and let $u_h\in\dgzero$ be
  its nonvariational finite element approximation. Assume that the
  numerical fluxes $\hat U$ and $\hat{\geovec p}$ are consistent then
  we have the following orthogonality result:
  \begin{equation}
    \bih{u_h-u}{\Psi} = J(\Psi) \Foreach \Psi\in\dgzero,
  \end{equation}
  with the error functional given by
  \begin{align}
    J(\Psi) = \int_\W \frob{\qp{\Hess u-\H[u]}}{\qp{\A\Psi}} \d \geovec x.
  \end{align}
\end{Lem}
\begin{Proof}
  Using the consistency result and that $\jump{u}=0$ we conclude
  \begin{align*}
    \bih{u_h-u}{\Psi}
    &= \bih{u_h}{\Psi} + \int_\W \frob{\A}{\H[u]} \Psi \d \geovec x
     = l(\Psi) + \int_\W \frob{\H[u]}{\qp{\A\Psi}} \d \geovec x\\
    &= -\int_\W \frob{\A}{\Hess u}\Psi - \frob{\H[u]}{\qp{\A\Psi}}\d \geovec x
     = J(\Psi),
  \end{align*}
  concluding the proof.
\end{Proof}
\begin{Rem}
If $\A$ is piecewise constant then since $\H[u]=\ltwoproj\qp{\Hess u}$ we
have $J(\Psi) = 0 $ and we recover the usual Galerkin orthogonality $\bih{u_h-u}{\Psi} = 0$.
We will show in the next section that in general the error functional $J$ is
of higher order for smooth enough $u$.
\end{Rem}


\section{Coercivity, continuity and convergence}
\label{sec:consistency}

In this section we examine the coercivity, continuity and
convergence of the method. 
We will focus on the fluxes given in Example~\ref{ex:ip-hessian} to
simplify the presentation. Furthermore we make the following additional assumption on
the problem data.
\begin{Hyp}[coercivity of the nonvariational problem]
  \label{hyp:coercivity-of-cont}
  For the rest of this section we will assume that the nonvariational
  operator is coercive, that is $\A\in\qb{\sob{k+1}{\infty}(\W)}^{d\times d}$ and that
  $\div\qp{\D\A} \leq 0$.
\end{Hyp}
\begin{Rem}[variational nature of the coercive problem]
  Under Assumption \ref{hyp:coercivity-of-cont} the
  problem can be written variationally. The solution to the
  nonvariational problem is the minimiser to the (degenerate) second
  order variational problem: Find $u\in\sobh{2}(\W)\cap\hoz(\W)$ such that
  \begin{equation}
    \cJ[u] = \inf_{v\in\sobh{2}(\W)\cap\hoz(\W)} \cJ[v],
  \end{equation}
  where 
  \begin{equation}
    \cJ[v] := \int_\W \qp{\frob{\A}{\Hess u} + f - \frob{\Hess}{\A}}u \d \geovec x.
  \end{equation}
\end{Rem}


\begin{Defn}[$\sobh{1}(\T{})$, $\sobh{2}(\T{})$ and $\sobh{-1}(\T{})$ norms]
  We introduce the broken $\sobh{1}(\T{})$ and $\sobh{2}(\T{})$ norms as
  \begin{gather}
    \enorm{u_h}{1}^2 
    := 
    \Norm{\nabla_h u_h}_{\leb{2}(\W)}^2
    +
    h^{-1} 
    \Norm{\jump{u_h}}_{\leb{2}(\E)}^2,
    \\
    \enorm{u_h}{2}^2 
    := 
    \Norm{\Hess_h u_h}_{\leb{2}(\W)}^2
    + 
    h^{-1} 
    \Norm{\jump{\nabla_h u_h}}_{\leb{2}(\E)}^2
    + 
    h^{-3} 
    \Norm{\jump{u_h}}_{\leb{2}(\E)}^2,
  \end{gather}
  and the $\sobh{-1}(\T{})$ norm as
  \begin{equation}
    \label{eq:neg-norm}
    \enorm{u_h}{-1} := \sup_{v_h\in\dgzero} \frac{\int_\W u_h v_h \d \geovec x}{\enorm{v_h}{1}}.
  \end{equation}
  These are equivalent to their continuous equivalent norms for
  functions in $\dg$.
\end{Defn}

\begin{Pro}[projection approximation in $\dg$]
  \label{pro:interpolation}
Let $\ltwoproj:\leb{2}(\W) \to \dg$ be the $\leb{2}(\W)$ orthogonal
projection operator defined by (\ref{eq:l2-proj-def}). Using standard
approximation arguments we have that
  \begin{equation}
    \begin{split}
      \enorm{v - \ltwoproj v}{1} &\leq C h^k\norm{v}_{\sobh{k+1}(\W)}\AND
      \\
      \Norm{v - \ltwoproj v}_{\leb{2}(\W)} &\leq C h^{k+1}\norm{v}_{\sobh{k+1}(\W)}.
    \end{split}
  \end{equation}
  In particular, let $\A_h$ denote the $\leb{2}$ orthogonal projection
  of $\A$ into the space of piecewise constant functions, then we have
  \begin{equation}
    \label{eq:approx-of-A-h}
    \Norm{\A - \A_h}_{\sob{1}{\infty}(\W)} 
    \leq
    \qp{1+C_1 h} \Norm{\D \A}_{\leb{\infty}(\W)}.
  \end{equation}
\end{Pro}

\begin{The}[stability of {$\H$} {\cite[Theorem 4.10]{Pryer:2012a}}]  
  \label{the:stability}
  Let $\H$ be defined as in Example \ref{ex:ip-hessian} then the
  dG Hessian is stable in the sense that
  \begin{equation}
    \label{eq:bound-the-liftings}
    \begin{split}
      \Norm{\Hess_h v_h - \H[v_h]}_{\leb{2}(\W)}^2
      &\leq
      C\qp{\int_\E h^{-1} 
      \norm{\jump{\nabla_h v_h}}^2
      +
      h^{-3}
      \norm{\jump{v_h}}^2 \d s
      }.
    \end{split}
  \end{equation}
  Consequently we have
  \begin{equation}
    \label{eq:stab-of-H}
    \Norm{\H[v_h]}_{\leb{2}(\W)}^2
    \leq
    C \enorm{v_h}{2}^2.
  \end{equation}
\end{The}

We now state the following technical Lemmata.

\begin{Lem}[upper bound on $\Norm{v_h\A - \ltwoproj\qp{v_h \A}}_{\leb{2}(\W)}$]
  \label{lem:neilan-tech}
  Let $\A\in\qb{\sob{k+1}{\infty}(\W)}^{d\times d}$ and $v_h \in \fes$ then it holds that
  \begin{gather}
    \Norm{\D_h\qp{v_h \A - \ltwoproj\qp{v_h \A}}}_{\leb{2}(\W)}
    \leq 
    C_2 h \qp{\qp{1+C_1h}\Norm{\D\A}_{\leb{\infty}(\W)} + \sum_{i = 2}^{k+1} \norm{\A}_{\sob{i}{\infty}(\W)}} \Norm{\nabla_h v_h}_{\leb{2}(\W)}.
  \end{gather}
\end{Lem}
\begin{Proof}
  Let $\A_h$ denote the $\leb{2}$ orthogonal projection of $\A$ into
  the space of piecewise constant functions. Then adding and
  subtracting appropriate terms we see
  \begin{equation}
    \Norm{\D_h\qp{v_h \A - \ltwoproj\qp{v_h \A}}}_{\leb{2}(\W)}
    =
    \Norm{\D_h\qp{v_h \qp{\A - \A_h} - \ltwoproj\qp{v_h {\A - \A_h}} }}_{\leb{2}(\W)}.
  \end{equation}
  Using the approximation
  properties of the $\leb{2}(\W)$ projection given in Proposition
  \ref{pro:interpolation} we see that
  \begin{equation}
    \begin{split}
      \Norm{\D_h\qp{v_h \A - \ltwoproj\qp{v_h \A}}}_{\leb{2}(\W)}
      &\leq Ch^{k} \norm{v_h \qp{\A -
          \A_h}}_{\sobh{k+1}(\W)}.
      \\
      &\leq 
      Ch^{k} \Norm{\A - \A_h}_{\sob{k+1}{\infty}(\W)}\Norm{v_h}_{\sobh{k+1}(\W)}.
    \end{split}
  \end{equation}
  Now using the properties of $\A_h$ (\ref{eq:approx-of-A-h}) and
  inverse inequalities we have
  \begin{equation}
    \begin{split}
      \Norm{\D_h\qp{v_h \A - \ltwoproj\qp{v_h \A}}}_{\leb{2}(\W)}
      &\leq Ch \Norm{\A - \A_h}_{\sob{k+1}{\infty}(\W)}\Norm{\nabla_h v_h}_{\leb{2}(\W)}
      \\
      &\leq 
      C h \qp{\qp{1+C_1h}\Norm{\D\A}_{\leb{\infty}(\W)} + \sum_{i = 2}^{k+1} \norm{\A}_{\sob{i}{\infty}(\W)}} \Norm{\nabla_h v_h}_{\leb{2}(\W)}
      ,
    \end{split}
  \end{equation}
  as required.
\end{Proof}

\begin{Lem}[upper bound on some skeletal terms]
  \label{lem:tech}
  Let $v_h\in\dg$, $\B\in\sob{1}{\infty}(\W)^{d\times d}$ and
  $\B_h\in\dg^{d\times d}$ be the $\leb{2}$ orthogonal projection of
  $\B$, then in view of trace and inverse inequalities we have the
  following:
  \begin{equation}
    \begin{split}
      \Norm{\avg{\D_h\ltwoproj\qp{v_h\B}}}_{\leb{2}(\E)} 
      &\leq
      C h^{-1/2} 
      \Norm{\D_h\qp{\ltwoproj\qp{v_h\B}}}_{\leb{2}(\W)}
      \\
      &\leq 
      C_3 h^{-1/2}
      \Norm{\B}_{\sob{1}{\infty}(\W)} \enorm{v_h}{1},
      \\
      \Norm{\jump{\ltwoproj\qp{v_h\B}}}_{\leb{2}(\E)}
      &\leq 
      C h^{1/2}
      \Norm{\D_h\qp{\ltwoproj\qp{v_h\B}}}_{\leb{2}(\W)}
      \\
      &\leq 
      C_4 h^{1/2}
      \Norm{\B}_{\sob{1}{\infty}(\W)} \enorm{v_h}{1},
      \\
      \Norm{\avg{\nabla_h v_h}}_{\leb{2}(\E)}
      &\leq
      C_5 h^{-1/2} \Norm{\nabla_h v_h}_{\leb{2}(\W)},
      \\
      \Norm{\avg{v_h}}_{\leb{2}(\E)}
      &\leq
      C_6 h^{1/2} \Norm{\nabla_h v_h}_{\leb{2}(\W)},
      \\
      \Norm{\jump{\B_h}}_{\leb{\infty}(\E)}
      &\leq
      C_7 h \Norm{\D\B}_{\leb{\infty}(\W)}.
    \end{split}
  \end{equation}
\end{Lem}
\begin{Proof}
  For brevity we prove only the first inequality, the second and third
  follow similar arguments. In view of the definition of the average
  operator (\ref{eqn:average}) it follows that
  \begin{equation}
    \Norm{\avg{\D_h\ltwoproj\qp{v_h\B}}}_{\leb{2}(\E)}
    \leq
    \frac{1}{2}
    \sum_{K\in\T{}}
    \Norm{\D_h \ltwoproj\qp{v_h\B}}_{\leb{2}(\partial K)}.
  \end{equation}
  Now by a trace inequality we see that
  \begin{equation}
    \Norm{\avg{\D_h\ltwoproj\qp{v_h\B}}}_{\leb{2}(\E)}
    \leq
    C
    \sum_{K\in\T{}}
    h^{-1/2} \Norm{\D_h \ltwoproj\qp{v_h\B}}_{\leb{2}(K)}.
  \end{equation}
  Using the stability of the $\leb{2}(\W)$ projection operator in
  $\sobh1(\W)$ we have
  \begin{equation}
    \begin{split}
      \Norm{\avg{\D_h\ltwoproj\qp{v_h\B}}}_{\leb{2}(\E)} &\leq C \sum_{K\in\T{}}
      h^{-1/2} \Norm{\D_h \qp{v_h\B}}_{\leb{2}(K)}
      \\
      &\leq
      C_3 h^{-1/2}
      \Norm{\B}_{\sob{1}{\infty}(\W)} \enorm{v_h}{1},
    \end{split}
  \end{equation}
  as required.
\end{Proof}

\begin{The}[discrete continuity and coercivity]
  \label{the:bound-and-coer}
  Let the conditions in Assumption \ref{hyp:coercivity-of-cont}
  hold. Suppose also that $\D\A$ is sufficiently small such that
  \begin{equation}
    \label{eq:restriction-on-A}
    \begin{split}
      \gamma 
      -
      2\epsilon 
      -
      \big(
      C_4C_5
      +
      C_1C_4 h
      +
      \frac{C_6C_7}{4}h
      +
      C_2h\qp{1+C_1h}
      \big)
      &\Norm{\D\A}_{\leb{\infty}(\W)}
      \\
      &-
      C_2 h
      \sum_{i=2}^{k+1} \norm{\A}_{\sob{i}{\infty}(\W)}
      > 0,
    \end{split}
    \end{equation}
    where $\gamma$ is the ellipticity constant, $C_1$ is a constant
    appearing in Proposition \ref{pro:interpolation} with $k=0$, $C_2,
    C_4 \AND C_5$ are the constants appearing in Lemmata
    \ref{lem:neilan-tech} and \ref{lem:tech} and $\epsilon >0$ is some
    parameter. In addition assume $\sigma$ the penalisation term is
    sufficiently large, specifically
  \begin{equation}
    \begin{split}
      \sigma -
      \frac{C_3^2\qp{1+C_1h}^2\Norm{\D\A}_{\leb{\infty}(\W)}^2 +
        4\epsilon^2}{4\epsilon} 
        &- \qp{\frac{C_4C_5\qp{1+C_1 h}}{2} + C_7h}\Norm{\D\A}_{\leb{\infty}(\W)} 
      \\ &\qquad - \frac{\qp{\theta + 1}^2
        C_5^2 \Norm{\avg{\A_h}}_{\leb{\infty}(\E)}^2} {4\epsilon} > 0.
    \end{split}
  \end{equation}
  where $\A_h$ is some piecewise constant approximation to $\A$.

  Then there exist positive constants $C_B \AND C_C$ such that
  \begin{gather}
    \norm{\bih{u_h}{v_h}} \leq C_B\enorm{u_h}{1}\enorm{v_h}{1} \AND
    \\
    \bih{u_h}{u_h} \geq C_C \enorm{u_h}{1}^2 \Foreach u_h, v_h \in\dgzero.
  \end{gather}
\end{The}

We postpone the proof of this theorem to the end of this section and first
prove the error estimates for our discrete solution.


The results of Theorem \ref{the:bound-and-coer} allow us to invoke
Strang's second Lemma.
\begin{Cor}[Strang {\cite[c.f.]{ErnGuermond:2004}}]
  \label{cor:strang}
  There exists a $C>0$ such that
  \begin{equation}
    \enorm{u - u_h}{1}
    \leq C
    \qp{
      \inf_{v_h\in\dgzero} \enorm{u-v_h}{1} 
      +
      \sup_{w_h\in\dgzero} \frac{\norm{\bih{u}{w_h} - l(w_h)}}{\enorm{w_h}{1}}
    }.
  \end{equation}
\end{Cor}

\begin{Lem}[discrete negative norm convergence]
  \label{lem:discrete-neg-norm}
  Let $\A\in\leb{\infty}(\W)$ and $u\in\sobh{k+3}(\W)$. Then we have
  that there exists a constant $C>0$ such that
  \begin{equation}
    \enorm{\frob{\A}{\qp{\Hess u - \H[u]}}}{-1} 
    \leq
    Ch^{k+1}\Norm{\A}_{\infty}\Norm{u}_{\sobh{k+3}(\W)} .
  \end{equation}
\end{Lem}
\begin{Proof}
  We have, in view of Cauchy--Schwarz inequality, that for
  $\Psi\in\dgzero$
  \begin{equation}
    \begin{split}
      \int_\W \frob{\A}{\qp{\Hess u - \H[u]}} \Psi \d \geovec x
      &\leq
      \Norm{\A}_{\leb{\infty}(\W)} \Norm{\Hess u - \H[u]}_{\leb{2}(\W)} \Norm{\Psi}_{\leb{2}(\W)}
      \\
      &\leq C h^{k+1} \Norm{\A}_{\leb{\infty}(\W)} \Norm{u}_{\sobh{k+3}(\W)} \enorm{\Psi}{1},
    \end{split}
  \end{equation}
  since $\H[u] = \ltwoproj \Hess u$ and by the definition of
  $\enorm{\cdot}{1}$. The result follows noting the definition of the
  discrete negative norm in (\ref{eq:neg-norm}).
\end{Proof}
\begin{Rem}
  Noting the definition of the error functional in the Galerkin
  orthogonality we deduce that
  \begin{equation}
    \enorm{J}{-1} = \Oh(h^{k+1}).
  \end{equation}
\end{Rem}

\begin{The}[convergence of the nonvariational method]
  Let $u$ solve the nonvariational problem
  (\ref{eq:nonvariational-problem}) and $\qp{u_h, \H[u_h]}$ solve the
  nonvariational finite element approximation
  (\ref{eq:nonvariational-fe-approx}) where $\H[u_h]$ is a consistent
  approximation of $\Hess u$ (for example that given in Example
  \ref{ex:ip-hessian}). Then the following error bound holds:
  \begin{gather}
    \label{eq:enorm-error-bound}
    \enorm{u - u_h}{1} \leq C \qp{ h^k \norm{u}_{\sobh{k+1}(\W)} + h^{k+1} \norm{u}_{\sobh{k+3}(\W)}}.
  \end{gather}
\end{The}
\begin{Proof}
  The proof of (\ref{eq:enorm-error-bound}) is immediate from applying
  Proposition \ref{pro:interpolation} to Corollary \ref{cor:strang}
  with 
  $v_h = \ltwoproj u$
  and noting that the bound for the
  consistency error nothing but the result of Lemma
  \ref{lem:discrete-neg-norm}, concluding the proof.
\end{Proof}


To conclude this section we prove Theorem \ref{the:bound-and-coer}.
\begin{Proof}[Theorem \ref{the:bound-and-coer}]
  Let $u_h, v_h \in\dg$, then we have 
  \begin{equation}
    \begin{split}
      \bih{u_h}{v_h} 
      &=
      -\int_\W \frob{\A}{\H[u_h]} v_h \d \geovec x
      +
      \sigma h^{-1} \int_{\E\cup\partial\W} \Transpose{\jump{u_h}}\jump{v_h}\d s
      \\
      &=
      -\int_\W \frob{\H[u_h]}{\qp{v_h\A}}\d \geovec x
      +
      \sigma h^{-1} \int_{\E\cup\partial\W} \Transpose{\jump{u_h}}\jump{v_h}\d s
      \\
      &=
      -\int_\W \frob{\H[u_h]}{\ltwoproj\qp{v_h\A}}\d \geovec x
      +
      \sigma h^{-1} \int_{\E\cup\partial\W} \Transpose{\jump{u_h}}\jump{v_h}\d s
      .
    \end{split}
  \end{equation}
  Now making use of the
  definition of $\H$ from Example \ref{ex:ip-hessian} we see
  \begin{equation}
    \begin{split}
      \bih{u_h}{v_h}
      &=
      -
      \int_\W 
      \frob{\Hess_h u_h}{\ltwoproj\qp{v_h\A}}\d \geovec x
      +
      \int_\E
      \frob{\tjump{\nabla_h u_h}}{\avg{\ltwoproj\qp{v_h\A}}}\d s
      \\
      &\qquad
      -
      \int_{\E\cup\partial\W}
      \theta \Transpose{\jump{u_h}}\avg{\D_h \qp{\ltwoproj\qp{v_h\A}}}\d s
      +
      \sigma h^{-1} \int_{\E\cup\partial\W} \Transpose{\jump{u_h}}\jump{v_h}\d s
      .
    \end{split}
  \end{equation}
  Adding and subtracting appropriate terms we have that
  \begin{equation}
    \begin{split}
      \bih{u_h}{v_h}
      &= 
      -
      \int_\W 
      \frob{\A}{\Hess_h u_h} v_h
      +
      \frob{\Hess_h u_h}{\qp{v_h\A - \ltwoproj\qp{v_h \A}}}\d \geovec x
      \\
      &\qquad
      +
      \int_\E
      \frob{\tjump{\nabla_h u_h}}{\avg{\ltwoproj\qp{v_h\A}}}\d s
      -
      \int_{\E\cup\partial\W}
      \theta\Transpose{\jump{u_h}}\avg{\D_h \qp{\ltwoproj\qp{v_h\A}}}\d s
      \\
      &\qquad +
      \sigma h^{-1} \int_{\E\cup\partial\W} \Transpose{\jump{u_h}}\jump{v_h}\d s
      ,
    \end{split}
  \end{equation}
  which rewriting variationally gives
  \begin{equation}
    \label{eq:bi-form}
    \begin{split}
      \bih{u_h}{v_h}
      &= 
      -
      \int_\W 
      \D_h\qp{\A \nabla_h u_h} v_h
      +
      \D_h \A \nabla_h u_h v_h
      +
      \frob{\Hess_h u_h}{\qp{v_h\A - \ltwoproj\qp{v_h \A}}}\d \geovec x
      \\
      &\qquad
      +
      \int_\E
      \frob{\tjump{\nabla_h u_h}}{\avg{\ltwoproj\qp{v_h\A}}}\d s
      -
      \int_{\E\cup\partial\W}
      \theta\Transpose{\jump{u_h}}\avg{\D_h \qp{\ltwoproj\qp{v_h\A}}}\d s
      \\
      &\qquad
      +
      \sigma h^{-1} \int_{\E\cup\partial\W} \Transpose{\jump{u_h}}\jump{v_h}\d s
      .
    \end{split}
  \end{equation}
  Note that 
  \begin{equation}
    \label{eq:var-eq}
    \begin{split}
      -
      \int_\W 
      &\D_h\qp{\A \nabla_h u_h} v_h
      +
      \frob{\Hess_h u_h}{\qp{v_h\A - \ltwoproj\qp{v_h \A}}} \d \geovec x
      \\
      &=
      \sum_{K\in\T{}}
      \bigg[
      \int_K \qp{\A \nabla_h u_h} \nabla_h v_h
      - \D_h\qp{v_h\A - \ltwoproj\qp{v_h\A}}\nabla_h u_h \d \geovec x
      \\
      &\qquad +
      \int_{\partial K} -\qp{\A \nabla_h u_h} v_h\geovec n
      +
      \qp{\qp{v_h\A - \ltwoproj\qp{v_h \A}}\nabla_h u_h }\geovec n \d s
      \bigg]
      \\
      &=
      \sum_{K\in\T{}}
      \bigg[
      \int_K \qp{\A \nabla_h u_h} \nabla_h v_h
      - \D_h\qp{v_h\A - \ltwoproj\qp{v_h\A}}\nabla_h u_h \d \geovec x
      \\
      &\qquad -
      \int_{\partial K}\qp{\ltwoproj\qp{v_h \A}\nabla_h u_h }\geovec n \d s
      \bigg]
      \\
      &=
      \int_\W \qp{\A \nabla_h u_h} \nabla_h v_h
      - \D_h\qp{v_h\A - \ltwoproj\qp{v_h\A}}\nabla_h u_h \d \geovec x
      \\
      &\qquad 
      -
      \int_{\E\cup\partial\W}
      \Transpose {\jump{\ltwoproj\qp{v_h\A}}}\avg{\nabla_h u_h} \d s
      -
      \int_{\E}
      \frob{\tjump{\nabla_h u_h}}{\avg{\ltwoproj\qp{v_h\A}}} \d s,
    \end{split}
  \end{equation}
  and hence we see that upon substituting (\ref{eq:var-eq}) into
  (\ref{eq:bi-form}) that
  \begin{equation}
    \label{eq:bi-final-form}
    \begin{split}
      \bih{u_h}{v_h} 
      &= 
      \int_\W \qp{\A \nabla_h u_h} \nabla_h v_h + \D \A
      \nabla_h u_h v_h - {\D_h\qp{v_h\A - \ltwoproj\qp{v_h \A}}} \nabla_h u_h \d \geovec x
      \\
      & \qquad - \int_{\E\cup\partial\W} \theta\Transpose{\jump{u_h}}\avg{\D_h
        \qp{\ltwoproj\qp{v_h\A}}}  +
      \Transpose{\jump{\ltwoproj\qp{v_h\A}}}\avg{\nabla_h u_h} \d s
      \\
      & \qquad
      +
      \sigma h^{-1} \int_{\E\cup\partial\W} \Transpose{\jump{u_h}}\jump{v_h} \d s.
    \end{split}
  \end{equation}
  We proceed by applying
  Cauchy--Schwartz componentwise to (\ref{eq:bi-final-form}) and estimating
  $\theta$ by $1$ arriving at
  \begin{equation}
    \begin{split}
      \bih{u_h}{v_h} 
      &\leq \Norm{\nabla_h
        u_h}_{\leb{2}(\W)}\qp{\Norm{\A}_{\leb{\infty}(\W)} \Norm{\nabla_h v_h}_{\leb{2}(\W)} +
      \Norm{\D\A}_{\leb{\infty}(\W)}\Norm{v_h}_{\leb{2}(\W)}}
      \\
      &
      \qquad 
      +
      \Norm{\D_h\qp{v_h\A - \ltwoproj\qp{v_h\A}}}_{\leb{2}(\W)}\Norm{\nabla_h u_h}_{\leb{2}(\W)}
      \\
      &\qquad +
      \Norm{\jump{u_h}}_{\leb{2}(\E)} \Norm{\avg{\D_h\qp{\ltwoproj\qp{v_h\A}}}}_{\leb{2}(\E)}
      \\
      &\qquad +
      \Norm{\jump{\ltwoproj{v_h\A}}}_{\leb{2}(\E)} \Norm{\avg{\nabla_h u_h}}_{\leb{2}(\E)}
      +
      \sigma h^{-1}\Norm{\jump{u_h}}_{\leb{2}(\E)} \Norm{\jump{v_h}}_{\leb{2}(\E)}.
    \end{split}
  \end{equation}
  In view of Lemma \ref{lem:neilan-tech} and the Poincar\'e inequality
  we have
  \begin{equation}
    \label{eq:bilinear-boundedness}
    \begin{split}
      \bih{u_h}{v_h} 
      &\leq
      \qp{
        \Norm{\A}_{\leb{\infty}(\W)}       
        +
        C_P \Norm{\D\A}_{\leb{\infty}(\W)}
        +
        C_2 h \Norm{\A}_{\sob{k+1}{\infty}(\W)} 
      }
      \Norm{\nabla_h u_h}_{\leb{2}(\W)}\Norm{\nabla_h v_h}_{\leb{2}(\W)}
      \\
      &
      \qquad 
      +\Norm{\jump{u_h}}_{\leb{2}(\E)} \Norm{\avg{\D_h\qp{\ltwoproj\qp{v_h\A}}}}_{\leb{2}(\E)}
      \\
      &\qquad +
      \Norm{\jump{\ltwoproj{v_h\A}}}_{\leb{2}(\E)} \Norm{\avg{\nabla_h u_h}}_{\leb{2}(\E)}
      +
      \sigma h^{-1}\Norm{\jump{u_h}}_{\leb{2}(\E)} \Norm{\jump{v_h}}_{\leb{2}(\E)}.
    \end{split}
  \end{equation} 
  For the skeletal terms we apply the result of Lemma \ref{lem:tech}
  which upon substituting these into (\ref{eq:bilinear-boundedness})
  we see that
  \begin{equation}
    \norm{\bih{u_h}{v_h}} \leq C_B \enorm{u_h}{1}\enorm{v_h}{1}
  \end{equation}
  as required.

  For coercivity we use the equality given in (\ref{eq:bi-final-form})
  with $v_h=u_h$ to find
  \begin{equation}
    \label{eq:bound-i-1-6}
    \begin{split}
      \bih{u_h}{u_h}
      &=
      \int_\W \qp{\A \nabla_h u_h} \nabla_h u_h + 
      \D \A
      \nabla_h u_h u_h - {\D_h\qp{u_h\A - \ltwoproj\qp{u_h \A}}} \nabla_h u_h \d \geovec x
      \\
      & \qquad - \int_{\E\cup\partial\W} 
      \theta\Transpose{\jump{u_h}}\avg{\D_h
        \qp{\ltwoproj\qp{u_h\A}}}  +
      \Transpose{\jump{\ltwoproj\qp{u_h\A}}}\avg{\nabla_h u_h} 
      \\
      &\qquad\qquad + \sigma h^{-1}\Transpose{\jump{u_h}}\jump{u_h} \d s
      \\
      &=
      \sum_{i=1}^6 \cI_i 
      .
    \end{split}
  \end{equation}
  We proceed by bounding each term individually. By the ellipticity of
  the problem we have that
  \begin{equation}
    \label{eq:bound-1}
      \cI_1 
      =
      \int_\W \qp{\A \nabla_h u_h} \nabla_h u_h  \d \geovec x
      \geq
      \gamma \Norm{\nabla_h u_h}_{\leb{2}(\W)}^2.
  \end{equation}
  By the coercivity of the problem we have
  \begin{equation}
    \label{eq:bound-2}
    \cI_2
    =
    \int_\W \D \A
    \nabla_h u_h u_h
    =
    \int_\W
    \D \A
    \frac{1}{2}\nabla_h \qp{u_h^2}
    \d \geovec x
    =
    -\frac{1}{2} \int_\W \div\qp{\D \A} u_h^2 \d \vec x > 0.
  \end{equation}
  By the Cauchy--Schwartz inequality and making use of Lemma \ref{lem:neilan-tech}
  \begin{equation}
    \begin{split}
      \label{eq:bound-3}
      - \cI_3 
      &=
      \int_\W{\D_h\qp{u_h\A - \ltwoproj\qp{u_h \A}}}
      \nabla_h u_h \d \geovec x
      \\
      &\leq
      \Norm{\D_h\qp{u_h\A - \ltwoproj\qp{u_h \A}}}_{\leb{2}(\W)} \Norm{\nabla_h u_h}_{\leb{2}(\W)}
      \\
      &\leq 
      C_2 h \qp{\qp{1+C_1h}\Norm{\D\A}_{\leb{\infty}(\W)} + \sum_{i = 2}^{k+1} \norm{\A}_{\sob{i}{\infty}(\W)}} \Norm{\nabla_h v_h}_{\leb{2}(\W)}^2.
    \end{split}
  \end{equation}
  We combine the fourth and fifth terms and let $\A_h$ denote the
  $\leb{2}$ orthogonal projection of $\A$ onto the space of piecewise
  constant functions. Upon adding and subtracting appropriate terms
  \begin{equation}
    \begin{split}
      - \cI_4 - \cI_5
      &=
      \int_{\E\cup\partial\W} \theta \Transpose{\jump{u_h}}\avg{\D_h
        \qp{\ltwoproj\qp{u_h\A}}}
      +
      \Transpose{\jump{\ltwoproj\qp{u_h\A}}}\avg{\nabla_h u_h} \d s
      \\
      &=
      \int_{\E\cup\partial\W} 
      \theta \qp{\Transpose{\jump{u_h}}\avg{\D_h
          \qp{\ltwoproj\qp{u_h\A - u_h\A_h}}}
        + 
        \theta\Transpose{\jump{u_h}}
        \avg{\D_h
          \qp{{u_h\A_h}}}
      }
      \\
      &\qquad+
      \Transpose{\jump{\ltwoproj\qp{u_h\A - u_h\A_h}}}\avg{\nabla_h u_h} 
      + 
      \Transpose{\jump{u_h\A_h}}\avg{\nabla_h u_h}
      \d s.
    \end{split}
  \end{equation}
  Using the identities
  \begin{gather}
    \int_{\E\cup\partial\W} \jump{u_h \A_h}\d s
    =
    \int_\E \avg{\A_h} \jump{u_h} \d s
    +
    \int_{\E\cup\partial\W} \Transpose{\jump{\A_h}} \avg{u_h} \d s \quad \AND
    \\
    \int_{\E\cup\partial\W} \avg{\A_h \nabla_h u_h} \d s
    =
    \int_\E \frac{1}{4}{\jump{A_h}}\jump{\nabla_h u_h} \d s
    +
    \int_{\E\cup\partial\W} \avg{\A_h}\avg{\nabla_h u_h} \d s
  \end{gather}
  we have that
  \begin{equation}
    \begin{split}
      - \cI_4 - \cI_5
      &=
      \int_{\E\cup\partial\W} 
      \theta \Transpose{\jump{u_h}}\avg{\D_h
          \qp{\ltwoproj\qp{u_h\A - u_h\A_h}}}
        \\
        &\qquad + 
        \Transpose{\jump{\ltwoproj\qp{u_h\A - u_h\A_h}}}\avg{\nabla_h u_h} 
      +
      \qp{\theta + 1}
      \qp{\avg{\A_h}{\jump{u_h}}}\avg{\nabla_h u_h}
        \\
      &\qquad
      +
      \jump{A_h} \avg{\nabla_h u_h} \avg{u_h}
      +
      \frac{\theta}{4}
      \Transpose{\jump{A_h}}\jump{u_h}\jump{\nabla_h u_h}
      \d s.    
    \end{split}
  \end{equation}
  Using Cauchy--Schwartz we see
  \begin{equation}
    \begin{split}
      -\cI_4-\cI_5 
      &\leq
      \Norm{\jump{u_h}}_{\leb{2}(\E)} 
      \Norm{\avg{\D_h \qp{\ltwoproj\qp{u_h\A - u_h\A_h}}}}_{\leb{2}(\E)} 
      \\
      &\qquad
      +
      \Norm{\jump{\ltwoproj\qp{u_h\A - u_h\A_h}}}_{\leb{2}(\E)} 
      \Norm{\avg{\nabla_h u_h}}_{\leb{2}(\E)} 
      \\
      &\qquad
      +
      \qp{\theta + 1}
      \Norm{\avg{\A_h}}_{\leb{\infty}(\E)} 
      \Norm{\jump{u_h}}_{\leb{2}(\E)} 
      \Norm{\avg{\nabla_h u_h}}_{\leb{2}(\E)} 
      \\
      &\qquad
      +
      \Norm{\jump{A_h}}_{\leb{\infty}(\E)} 
      \Norm{\avg{\nabla_h u_h}}_{\leb{2}(\E)} 
      \Norm{\avg{u_h}}_{\leb{2}(\E)} 
      \\
      &\qquad
      +
      \frac{1}{4}
      \Norm{\jump{A_h}}_{\leb{\infty}(\E)} 
      \Norm{\jump{u_h}}_{\leb{2}(\E)} 
      \Norm{\jump{\nabla_h u_h}}_{\leb{2}(\E)}.
    \end{split}
  \end{equation}
  Making use of the various bounds from Lemma \ref{lem:tech} we have
  \begin{equation}
    \begin{split}
      -\cI_4-\cI_5 
      &\leq
      C_3h^{-1/2}\Norm{A - A_h}_{\sob{1}{\infty}(\W)}
      \enorm{u_h}{1}
      \Norm{\jump{u_h}}_{\leb{2}(\E)}
      \\
      &\qquad 
      +
      C_4C_5 \Norm{A - A_h}_{\sob{1}{\infty}(\W)}
      \enorm{u_h}{1}
      \Norm{\nabla_h u_h}_{\leb{2}(\W)}
      \\
      &\qquad 
      +
      \qp{\theta + 1}
      C_5 h^{-1/2}
      \Norm{\avg{\A_h}}_{\leb{\infty}(\E)} 
      \Norm{\jump{u_h}}_{\leb{2}(\E)} 
      \Norm{{\nabla_h u_h}}_{\leb{2}(\W)} 
      \\
      &\qquad 
      +
      \frac{C_5C_6C_7 h}{4} 
      \Norm{\D\A}_{\leb{\infty}(\W)}
      \Norm{{\nabla_h u_h}}_{\leb{2}(\W)}^2
      +
      C_7 \Norm{\D\A}_{\leb{\infty}(\W)}
      \Norm{\jump{u_h}}_{\leb{2}(\E)}^2
      \\
      &
      =:
      \sum_{i=1}^5 \cK_i
      .
    \end{split}
  \end{equation}
  We now apply a Cauchy inequality and use the approximation
  properties of $\A_h$ from Proposition \ref{pro:interpolation} to
  find for any $\epsilon > 0$ that
  \begin{equation}
    \label{eq:bound-k-1}
    \begin{split}
      \cK_1
      &=
      C_3h^{-1/2}\Norm{A - A_h}_{\sob{1}{\infty}(\W)}
      \enorm{u_h}{1} \Norm{\jump{u_h}}_{\leb{2}(\E)} 
      \\
      &\leq
      \frac{C_3^2\Norm{A - A_h}_{\sob{1}{\infty}(\W)}^2}{4\epsilon}
      h^{-1}\Norm{\jump{u_h}}_{\leb{2}(\E)}^2
      +
      \epsilon 
      \enorm{u_h}{1}^2
      \\
      &
      \leq
      \frac{C_3^2\Norm{A - A_h}_{\sob{1}{\infty}(\W)}^2 + 4\epsilon^2}{4\epsilon}
      h^{-1}\Norm{\jump{u_h}}_{\leb{2}(\E)}^2
      +
      \epsilon 
      \Norm{\nabla_h u_h}_{\leb{2}(\W)}^2
      \\
      &
      \leq
      \frac{C_3^2\qp{1+C_1h}^2\Norm{\D\A}_{\leb{\infty}(\W)}^2 + 4\epsilon^2}{4\epsilon}
      h^{-1}\Norm{\jump{u_h}}_{\leb{2}(\E)}^2
      +
      \epsilon 
      \Norm{\nabla_h u_h}_{\leb{2}(\W)}^2.
    \end{split}
  \end{equation}
  The other terms are bounded similarly in that
  \begin{equation}
    \begin{split}
      \cK_2 
      &=
      C_4C_5 \Norm{A - A_h}_{\sob{1}{\infty}(\W)}
      \enorm{u_h}{1}
      \Norm{\nabla_h u_h}_{\leb{2}(\W)}
      \\
      &\leq 
      \frac{C_4C_5 \Norm{A - A_h}_{\sob{1}{\infty}(\W)}}{2}
      \qp{
        \enorm{u_h}{1}^2
        +
        \Norm{\nabla_h u_h}_{\leb{2}(\W)}^2
      }
      \\
      &\leq
      C_4C_5\qp{1+C_1 h}\Norm{\D\A}_{\leb{\infty}(\W)}
      \Norm{\nabla_h u_h}_{\leb{2}(\W)}^2
      \\
      &\qquad 
      +
      \frac{C_4C_5\qp{1+C_1 h}\Norm{\D\A}_{\leb{\infty}(\W)}}{2}
      h^{-1}\Norm{\jump{u_h}}_{\leb{2}(\E)}^2,
    \end{split}
  \end{equation}
  and
  \begin{equation}
    \begin{split}
      \cK_3
      &=
      \qp{\theta + 1}
      C_5 h^{-1/2}
      \Norm{\avg{\A_h}}_{\leb{\infty}(\E)} 
      \Norm{\jump{u_h}}_{\leb{2}(\E)} 
      \Norm{{\nabla_h u_h}}_{\leb{2}(\W)} 
      \\
      &\leq
      \epsilon \Norm{{\nabla_h u_h}}_{\leb{2}(\W)}^2
      +
      \frac{\qp{\theta + 1}^2
        C_5^2
        \Norm{\avg{\A_h}}_{\leb{\infty}(\E)}^2}
      {4\epsilon}
      h^{-1}
      \Norm{\jump{u_h}}_{\leb{2}(\E)}^2.
    \end{split}
  \end{equation}
  Note that the final two terms are already in their desired form since
  \begin{gather}
    \cK_4
    =
    \frac{C_5C_6C_7 h}{4} 
    \Norm{\D\A}_{\leb{\infty}(\W)}
    \Norm{{\nabla_h u_h}}_{\leb{2}(\W)}^2
    \\
    \label{eq:bound-k-5}
    \cK_5 = 
    C_7 \Norm{\D\A}_{\leb{\infty}(\W)}
    \Norm{\jump{u_h}}_{\leb{2}(\E)}^2.
  \end{gather}
  Collecting the bounds from
  (\ref{eq:bound-k-1})--(\ref{eq:bound-k-5}) shows
  \begin{equation}
    \label{eq:bound-4}
    \begin{split}
      -\cI_4-\cI_5 
      &\leq
      \qp{2\epsilon + C_5\qp{C_4 + h\qp{C_1C_4 + \frac{C_6C_7}{4}}}
      \Norm{\D\A}_{\leb{\infty}(\W)}}
      \Norm{{\nabla_h u_h}}_{\leb{2}(\W)}^2
      \\
      &
      \qquad 
      +
      \Bigg(
        \frac{C_3^2\qp{1+C_1h}^2\Norm{\D\A}_{\leb{\infty}(\W)}^2 + 4\epsilon^2}{4\epsilon}
        \\
        &\qquad\qquad\qquad +
        \qp{\frac{C_4C_5\qp{1+C_1 h}}{2} + C_7h}\Norm{\D\A}_{\leb{\infty}(\W)}
        \\
        &\qquad\qquad\qquad\qquad\qquad +
        \frac{\qp{\theta + 1}^2
          C_5^2
          \Norm{\avg{\A_h}}_{\leb{\infty}(\E)}^2}
        {4\epsilon}
      \Bigg)
      h^{-1}
      \Norm{\jump{u_h}}_{\leb{2}(\E)}^2.
    \end{split}
  \end{equation}
  The final term in (\ref{eq:bound-i-1-6}) is given by
  \begin{equation}
    \label{eq:bound-6}
    \cI_6 
    =
    \int_{\E\cup\partial\W} \sigma h^{-1}\Transpose{\jump{u_h}}\jump{u_h} \d s 
    =
    \sigma h^{-1} \Norm{\jump{u_h}}_{\leb{2}(\E)}^2.
  \end{equation}
  Finally, collecting the bounds from (\ref{eq:bound-1}),
  (\ref{eq:bound-2}), (\ref{eq:bound-3}), (\ref{eq:bound-4}) and
  (\ref{eq:bound-6}) shows

  \begin{equation}
    \label{eq:coercivity-proof}
    \begin{split}
      \bih{u_h}{u_h} 
      &\geq 
      \Bigg( \gamma 
      -
      2\epsilon 
      -
      C_2h
      \qp{
        \qp{1+C_1h}\Norm{\D\A}_{\leb{\infty}(\W)} + \sum_{i=2}^{k+1} \norm{\A}_{\sob{i}{\infty}(\W)}
      }
      \\
      &\qquad\qquad -
      C_5\qp{C_4 +h}\qp{C_1C_4+\frac{C_6C_7}{4}}
      \Norm{\D\A}_{\leb{\infty}(\W)}
      \Bigg)
    \Norm{\nabla_h u_h}_{\leb{2}(\W)}^2 
    \\
    &\qquad +
    \Bigg(
    \sigma - 
    \frac{C_3^2\qp{1+C_1h}^2\Norm{\D\A}_{\leb{\infty}(\W)}^2 + 4\epsilon^2}{4\epsilon}
    \\
    &\qquad\qquad\qquad -
    \qp{\frac{C_4C_5\qp{1+C_1 h}}{2} + C_7h}\Norm{\D\A}_{\leb{\infty}(\W)}
    \\
    &\qquad\qquad\qquad\qquad\qquad -
    \frac{\qp{\theta + 1}^2
      C_5^2
      \Norm{\avg{\A_h}}_{\leb{\infty}(\E)}^2}
    {4\epsilon}
    \Bigg)
    h^{-1}
    \Norm{\jump{u_h}}_{\leb{2}(\E)}^2.
  \end{split}
\end{equation}
Coercivity of the discrete bilinear form follows using the assumption
in Theorem \ref{the:bound-and-coer}, by choosing $\epsilon$
sufficiently small and the penalisation parameter $\sigma$
sufficiently large for small enough $h$.
\end{Proof}

\begin{Rem}[the coercivity bound]
  We note that the coercivity bound relies on the term $\D\A$ not
  becoming too large, as specified in Theorem \ref{the:bound-and-coer}. If
  it is we view this as an advection dominated problem. Our
  numerical experiments suggest that there is \emph{no} condition on the
  size of this term.
  
  If the coefficient matrix $\A$ is divergence free, \ie $\D\A =
  \geovec 0$ then the bound simplifies considerably. For example, in
  the case that $\A$ is constant we regain the same theoretical
  results as for the method given in Example
  \ref{eq:laplace-formulation}.
\end{Rem}


\section{Numerical experiments}
\label{sec:numerics}

In this section we detail numerical experiments carried out in the
finite element package \dunefem 
\cite{DednerKloefkornNolteOhlberger:2010}
which is based on the \dune software framework
\cite{BastianBlattDednerEngwerKlofkornKornhuberOhlbergerSander:2008,
  BastianBlattDednerEngwerKlofkornOhlbergerSander:2008}. 
The code used to test the method will be made freely available within the
\dunefem-Howto in a future release.

We present some benchmark problems designed such that the exact
solution is known. In each of the experiments the domain $\W =
[0,1]^2$ and we consider the coefficient matrix to be 
\begin{equation}
  \label{eqn:benchmark-model-operator}
  \A(\vec x) =
  \begin{bmatrix}
    1 &  b(\vec x)
    \\
    b(\vec x) & a(\vec x)
  \end{bmatrix}
\end{equation}
varying $a(\vec x)$ and $b(\vec x)$. 

In each of the numerical experiments we make use a stabilised
conjugate gradient solver taken from the \duneistl module
\cite{BlattBastian:2007}
preconditioned with an incomplete LU
factorisation. We choose the penalty parameter $\sigma = 20$.

\subsection{Test 1 : a coercive operator} 
\label{test:1-coercive}

In this test we take the components of $\A$ such that the operator is
coercive, fitting into the analytical framework presented in
\S\ref{sec:consistency}. With $\geovec x = \qp{x_1,x_2}$, we set
\begin{align}
  a(\geovec x) &= -\ln\qp{\qp{x_1-1/2}^2+10^{{-10}}}+1
  \\
  b(\geovec x) &= 0.
\end{align}

We choose the problem data such that the exact solution is given by 
\begin{equation}
  u(\geovec x) = \sin{\pi x_1} \sin{\pi x_2}
\end{equation}
and approximate this using the formulation
(\ref{eq:nonvariational-fe-approx}). In Tables
\ref{subtbl:test1-k-1}--\ref{subtbl:test1-k-2} we present
the results for the cases $k=1,2$, numerically demonstrating that the
analytical rates of convergence are achieved in the dG energy norm,
moreover, optimal convergence is achieved in $\leb{2}(\W)$.

\begin{table}[h]
  \centering
  \caption{  \label{tbl:test1-coercive}
    \ref{test:1-coercive} - Test 1. We present errors and convergence rates of the approximation given by solving (\ref{eq:nonvariational-fe-approx}).}
  \subfloat[\label{subtbl:test1-k-1} Piecewise linears, $k=1$.]{
    \centering
      \small
      \begin{tabular}{c|c|c|c|c}
        \hline 
        $\# \text{elements}$ & $\Norm{u-u_h}$ & EOC & $\enorm{u-u_h}{1}$ & EOC \\ \hline 
        128 &   0.0196123 & 1.86116 &  0.414643& 0.953491   \\
        512 &   0.00506166 & 1.95408 &  0.209225& 0.986817   \\
        2048 &   0.00128044 & 1.98298 &  0.104907& 0.995937   \\
        8192 &  0.000321803 & 1.99238 &  0.0525047& 0.998597   \\
        32768 &  8.06862e-05 & 1.99578 &  0.0262623 & 0.999456   
      \end{tabular}
    }
    \\
    \subfloat[\label{subtbl:test1-k-2} Piecewise quadratics, $k=2$.]{
    \centering
      \small
      \begin{tabular}{c|c|c|c|c}
        \hline 
        $\# \text{elements}$ & $\Norm{u-u_h}$ & EOC & $\enorm{u-u_h}{1}$ & EOC \\ \hline 
        128 & 0.000475513 & 2.96865 &  0.0308463 &1.95   \\
        512 &  5.99935e-05 &2.98661 &  0.00779373 &1.98471 \\  
        2048 &   7.52887e-06& 2.9943 &  0.00195531 &1.99492  \\ 
        8192 &  9.42737e-07 &2.99751  & 0.000489487 &1.99805   \\
        32768 &  1.17929e-07& 2.99893  & 0.000122443 &1.99916   
      \end{tabular}
    }
\end{table}

\subsection{Test 2 : nondifferentiable operator \cite[\S 4.4]{LakkisPryer:2011}}  
\label{test:2-advection}

In this test we take $\A$ such that it is comparible to \cite[\S
4.4]{LakkisPryer:2011}. We take
\begin{align}
  a(\geovec x) &= 2
  \\
  b(\geovec x) &= \qp{x_1^2x_2^2}^{1/3}.
\end{align}

We choose the exact solution as in \ref{test:1-coercive} and conduct
the same tests. Tables \ref{subtbl:test2-k-1}--\ref{subtbl:test2-k-2}
detail the results. Note that this is not a coercive operator and as
such, does not fit into the analytical framework presented in
\S\ref{sec:consistency}, we do however still achieve optimal
convergence in $\Norm{\cdot}$ and $\enorm{\cdot}{1}$.

\begin{table}[h]
  \centering
  \caption{  \label{tbl:test2-advection}
    \ref{test:2-advection} - Test 2. We present errors and convergence rates of the approximation given by solving (\ref{eq:nonvariational-fe-approx}).}
  \subfloat[\label{subtbl:test2-k-1} Piecewise linears, $k=1$.]{
    \centering
      \small
      \begin{tabular}{c|c|c|c|c}
        \hline 
        $\# \text{elements}$ & $\Norm{u-u_h}$ & EOC & $\enorm{u-u_h}{1}$ & EOC \\ \hline 
        128 & 0.0172648 &1.89433 & 0.41799 &0.955709   \\
        512 & 0.00441656 &1.96684 & 0.210818 &0.987469   \\
        2048 & 0.00111269 &1.98887 & 0.105688 &0.996186   \\
        8192 & 0.000278969 &1.99588 & 0.0528915 &0.998707   \\
        32768 & 6.98234e-05 &1.99832 & 0.0264548 &0.999507   
      \end{tabular}
    }
    \\
    \subfloat[\label{subtbl:test2-k-2} Piecewise quadratics, $k=2$.]{
    \centering
      \small
      \begin{tabular}{c|c|c|c|c}
        \hline 
        $\# \text{elements}$ & $\Norm{u-u_h}$ & EOC & $\enorm{u-u_h}{1}$ & EOC \\ \hline 
        128 &  0.00047216 &2.9534   &0.0309416 &1.9514 \\  
        512 &  5.98325e-05& 2.98028  & 0.00782197 &1.98394 \\  
        2048 &  7.52118e-06& 2.9919   &0.00196325 &1.99429   \\
        8192 &  9.42426e-07 &2.99651   &0.000491575& 1.99776  \\
        32768 & 1.17933e-07 &2.99841  & 0.000122975& 1.99904   
      \end{tabular}
    }
\end{table}

\subsection{Test 3 : irregular solutions}
\label{test:3-irregular}

In this test we consider the case the exact solution does not satisfy
the regularity requirements presented in the analytical framework of
\S\ref{sec:consistency}, \ie $u\not\in\sobh{k+3}(\W)$. In addition we
consider the case that $u\not\in\sobh{2}(\W)$, demonstrating the
method converges even for viscocity solutions of the problem.

We consider the coercive operator from \S\ref{test:1-coercive} and
choose the problem data such that
\begin{equation}
  \label{eq:test3-solution-not-h3}
  u(\geovec x) 
  =
  \begin{cases}
    \frac{1}{4}\Big(\cos{8\pi\norm{\geovec x-\frac{1}{2}}^2}+1\Big) & \text{ if } \norm{\geovec x-\frac{1}{2}}^2 \leq \frac{1}{8}
    \\
    0 & \text{ otherwise }.
  \end{cases}
\end{equation}
Note that this function is $\sobh2(\W)$ but not $\sobh3(\W)$. We also
take the problem data such that
\begin{equation}
  \label{eq:test3-solution-not-h2}
  u(\geovec x)
  =
  \frac{100x_1(1-x_1)x_2(1-x_2)}{\norm{\geovec x}}.
\end{equation}
This function is $\sobh1(\W)$ but not $\sobh2(\W)$. The results are
given in Tables \ref{subtbl:test3-noth3}--\ref{subtbl:test3-noth2}.

In the case $u$ is given by (\ref{eq:test3-solution-not-h3}) the
scheme converges with optimal rate in the $\enorm{\cdot}{1}$ norm even
if the solution is not in $\sobh{3}$. The convergence in the $\leb{2}$ is more
erratic, but we observe the same behavior testing the standard IP FEM
taking $\A$ to be the identity.

In the case $u$ is given by (\ref{eq:test3-solution-not-h2}) the
convergence rates are suboptimal since the solution is not $\sobh{2}$.

\begin{table}[h]
  \centering
  \caption{  \label{tbl:test3-irregular}
    \ref{test:3-irregular} - Test 3. We present errors and convergence rates of the approximation given by solving (\ref{eq:nonvariational-fe-approx}). In both cases we consider $k=1$.}
  \subfloat[\label{subtbl:test3-noth3} The solution here is given in (\ref{eq:test3-solution-not-h3}). The function $u\in\sobh2(\W)$ but $u\not\in\sobh3(\W)$.]{
    \centering
      \small
      \begin{tabular}{c|c|c|c|c}
        \hline 
        $\# \text{elements}$ & $\Norm{u-u_h}$ & EOC & $\enorm{u-u_h}{1}$ & EOC \\ \hline 
        128 & 0.0362651 & 2.47943 &  0.837082& 0.939619\\
        512 & 0.0267684 & 0.43805  & 0.406003& 1.04388   \\
        2048 & 0.0179914 & 0.573227 &  0.253977& 0.67679   \\
        8192 & 0.00292357 & 2.6215  & 0.103168 &1.29971   \\
        32768 & 0.00174473 & 0.744729&   0.0541648 &0.929566\\   
        131072 & 0.000421749 & 2.04854&   0.0258935 &1.06476   
      \end{tabular}
    }
    \\
    \subfloat[\label{subtbl:test3-noth2} The solution here is given in (\ref{eq:test3-solution-not-h2}). The function $u\in\sobh1(\W)$ but $u\not\in\sobh2(\W)$.]{
    \centering
      \small
      \begin{tabular}{c|c|c|c|c}
        \hline 
        $\# \text{elements}$ & $\Norm{u-u_h}$ & EOC & $\enorm{u-u_h}{1}$ & EOC \\ \hline 
        128 & 0.223469 &1.80378 & 6.42181 &0.843123   \\
        512 & 0.0616572 &1.85773 &  3.49469& 0.877816   \\
        2048 & 0.017159 &1.84531  & 1.87984 &0.894556   \\
        8192 &  0.00509901& 1.75067&   1.00295& 0.906363  \\ 
        32768 & 0.00177874 &1.51936 &  0.531521& 0.916047   \\
        131072 &  0.00076433& 1.21859&   0.280092& 0.924224   
      \end{tabular}
    }
\end{table}

\section{Conclusions and outlook}

In this work we have extended the framework from
\cite{LakkisPryer:2011} for linear nonvariational problems to
incorporate discontinuous approximations.

We have shown the method presented (and subsequently that of the
continuous case from \cite{LakkisPryer:2011}) is well posed and
converges optimally under coercivity assumptions on the coefficient
matrix $\A$.

In the numerical experiments we note the the method is well posed and
converges optimally even for $\A$ which do not satisfy the coercivity
assumptions or $u$ which do not satisfy the regularity needed in the
analytical framework. This motivates that another analytical approach
needs to be developed. This approach can not be variational in nature
as such will be completely non standard. This is the topic of ongoing
research.


\bibliographystyle{alpha}
\bibliography{tristansbib,tristanswritings}

\end{document}